\documentclass[11pt]{amsart}
\usepackage{eucal}
\usepackage[all]{xy}
\usepackage{amsfonts,amssymb}
\usepackage{epsfig}
\usepackage[usenames]{color}
\usepackage{float} 
\xyoption{arc}

\newtheorem{theorem}{Theorem}[section]
\newtheorem{lemma}[theorem]{Lemma}

\newtheorem{corollary}[theorem]{Corollary}
\newtheorem{defn}[theorem]{Definition}

\newcommand{\Q}{\mathbb Q}
\newcommand{\C}{\mathbb C}
\newcommand{\R}{\mathbb R}
\newcommand{\K}{\mathbb K}

\newcommand{\Z}{\mathbb Z}

\newcommand{\bracket}[1]{ \{#1\} }

\title{Local Floer homology and infinitely many simple Reeb orbits}
\author{Mark McLean}

\begin{document}

\begin{abstract}
Let $Q$ be a Riemannian manifold such that the Betti numbers of its free loop space
with respect to some coefficient field are unbounded.
We show that every contact form on its unit contangent bundle supporting the natural contact structure
has infinitely many simple Reeb orbits.
This is an extension of a theorem by Gromoll and Meyer.
We also show that if a compact manifold admits a Stein fillable contact structure then there is a
possibly different such structure which also has infinitely many simple Reeb orbits for every supporting contact form.
We use local Floer homology along with symplectic homology to prove these facts.
\end{abstract}

\maketitle

\bibliographystyle{halpha}


\tableofcontents

\section{Introduction}

In this paper we are interested in contact manifolds that are the boundary $\partial M$
of certain symplectic manifolds $M$ called Liouville domains which will be defined later on.
We are interested in the Reeb orbits of such contact manifolds.
Let $\alpha_M$ be a supporting contact form on $\partial M$.
A {\it Reeb orbit} of period $T > 0$ is a smooth map $o : S^1 = \R / T \Z \rightarrow \partial M$ such that the vector
$\frac{d}{dt}(o(t))$ is in the kernal of $d\alpha_M$ and so that $\frac{d}{dt}(o(t))(\alpha_M)=1$.
Such an orbit is {\it simple} if the map $o$ is injective.
Suppose we choose a trivialization $\tau$ of the canonical bundle
of $M$ up to homotopy and a class $b \in H^2(M,\Z / 2\Z)$ then for every coefficient field $\K$
we can define a graded $\K$ vector space $SH_*(M,\K)$ called symplectic homology depending on $\tau$ and $b$.
The trivialization $\tau$ tells us how to grade the group $SH_*(M,\K)$.
The differential used to define symplectic homology involves counting
solutions to a certain differential equation called the perturbed Cauchy-Riemann equation.
Each solution has a sign $+$ or $-$ and we count the solutions with sign.
Different choices of the class $b$ will give different choices of sign,
and hence the differential changes when $b$ changes.
Symplectic homology is an invariant of $M$ up to deforming $M$ through Liouville domains
(assuming we do not change $(\tau,b)$).
There is also another invariant $\Gamma(M)$ called the {\it growth rate}.
The main theorem in this paper is the following:

\begin{theorem} \label{theorem:mainresult}
Suppose that $\partial M$
has only finitely many simple Reeb orbits, then:
\begin{enumerate}
\item There is a constant $C$ such that the rank of
$SH_k(M,\K)$ is bounded above by $C$ for all $k \notin [1-n,n]$
where $n$ is half the dimension of $M$.
\item $\Gamma(M) \leq 1$.
\end{enumerate}
This is true with respect to any coefficient field $\K$ and any choice of $\tau$ and $b$.
\end{theorem}

We do not require a genericty assumption here. These simple Reeb orbits
can be very degenerate. We prove this theorem by using work from
\cite{GinzburgGurel:localfloer}.
First of all we associate two Floer homology groups $CH_*(\gamma)$ and
$SH_*(\gamma)$ to each Reeb orbit.
Next we show that $SH_*(\gamma)$ is bounded above by $CH_*(\gamma) \oplus CH_{*-1}(\gamma)$.
We then put a bound on the rank of $CH_*(\gamma)$ by using work from \cite{GinzburgGurel:localfloer}.
It turns out that $SH_*(M,\K)$ is bounded above by the sum of $SH_*(\gamma)$
over all Reeb orbits plus the rank of $H^{n-*}(M)$ and this gives us our result.

For a compact Riemannian manifold $Q$, we define its {\it unit disk bundle}
$D^*Q$ to be the set of cotangent vectors of length $\leq 1$.
This is naturally a Liouville domain.
Its boundary $S^* Q$ is a contact manifold called the {\it unit cotangent bundle}.
By using the results from \cite{salamonweber:loop},
\cite{AbbondandoloSchwartz:cotangentloop}
or \cite{Viterbo:functorsandcomputations2} we get that
$SH_*(D^*Q,\K) = H_*(Q^{S^1},\K)$ where $Q^{S^1}$ is the free loop space of $Q$.
Hence we get the following corollary:
\begin{corollary} \label{corollary:gromollmeyergeneralization}
Suppose $H_k(Q^{S^1},\K)$ is unbounded for $k > n$ then
every contact form supporting the contact structure on $S^* Q$ has infinitely many simple Reeb orbits.
For instance if $Q$ is simply connected and its cohomology ring has at least two generators
over $\Q$ then $S^* Q$ has this property when $\K = \Q$
\cite{ViguePoirrieSullivan:closedgeodesic}.
\end{corollary}

This is an extension of a theorem by Gromoll and Meyer
\cite{GromollMeyer:periodicgeodesics}.
This corollary will also be proven in \cite{HryniewiczMacarini:localcontacthomology}
using similar methods. The main difference is that they use
contact homology which is the equivariant version of symplectic homology.
%
One can ask if other contact manifolds have infinitely many Reeb orbits.

\begin{theorem} \label{theorem:exoticliouvilledomain}

In each even dimension greater than $6$ there is a Liouville domain $M$
diffeomorphic to the ball such that $SH_*(M,\Q)$ has infinite rank in each degree.
\end{theorem}
We prove this in section \ref{section:exoticliouvilledomain}.
The {\it connected sum} $M \# N$ of two Liouville domains $M$ and $N$
is a new Liouville domain obtained by attaching a special $1$ handle called a Weinstein $1$-handle joining both.
By \cite{Cieliebak:handleattach}, we have $SH_*(M \# N,\K) = SH_*(M,\K) \oplus SH_*(N,\K)$.

By Theorem \ref{theorem:exoticliouvilledomain} and
Theorem \ref{theorem:mainresult}
we get the following corollary:
\begin{corollary} \label{corollary:reeborbitsoncontactstructures}
Let $N$ be any Liouville domain of dimension greater than $6$ with trivial first Chern class.
The boundary $\partial N$ admits a (possibly different) contact structure
with the property that every supporting contact form has infinitely many simple Reeb orbits.
This new contact structure is homotopic to the old one through hyperplane fields.
\end{corollary}
\proof of \ref{corollary:reeborbitsoncontactstructures}.
Let $M$ a Liouville domain as in Theorem \ref{theorem:exoticliouvilledomain} of the same dimension
as $N$
then $SH_*(M \# N,\Q)$ is infinitely generated in each degree.
Also $\partial(M \# N)$ is diffeomorphic to $\partial(N)$ and this diffeomorphism preserves the
homotopy type of the contact plane field within the space of hyperplane fields
(see \cite[Lemma 2.18]{McLean:computability}).
Hence by Theorem \ref{theorem:mainresult} we get that $\partial (M \# N)$ must have infinitely
many Reeb orbits for any supporting contact form.
\qed

\bigskip

{\bf Note:} While I was writing this paper I found out that
Leonardo Macarini and Umberto Hryniewicz were writing a similar
paper using contact homology instead of symplectic homology.
Even though the results are similar, I think that writing a version
of this paper from the perspective of symplectic homology is
interesting in its own right.

{\bf Acknowledgements:}
I would like to thank Viktor Ginzburg for his useful comments.
The author was partially supported by
NSF grant DMS-1005365.

\section{Definition of our Floer homology groups}

\subsection{Symplectic homology} \label{section:symplectichomology}

Symplectic homology in \cite{Viterbo:functorsandcomputations} was used to study Reeb orbits
on the boundary of Liouville domains. In this section we will define this invariant.
A {\it Liouville domain} is a compact manifold $M$ 
with boundary and a $1$-form $\theta_M$ satisfying:
\begin{enumerate}
\item $\omega_M := d\theta_M$ is a symplectic form.
\item The $\omega_M$-dual of $\theta_M$ is transverse
to $\partial M$ and pointing outwards.
\end{enumerate}
The boundary $\partial M$ is a contact manifold with contact form
$\alpha_M := \theta_M|_{\partial M}$.
Two Liouville domains are deformation equivalent if there is a smooth family
of Liouville domains joining them together.
Let $N$ be a Liouville domain with $c_1(N) = 0$.
We make some additional choices $\eta := (\tau,b)$ for $N$.
The element $\tau$ is a choice of trivialization of the canonical
bundle of $N$ up to homotopy and $b$ is an element of $H^2(N,\Z / 2\Z)$.
We will assume that $\partial N$ has discrete period spectrum ${\mathcal P}_N \subset \R$
(the set of periods of Reeb
orbits of $(\partial N,\alpha_{N})$).
For each pair of numbers $c<d$ where $c,d \in [-\infty,\infty]$
we will define a symplectic homology group $SH_*^{(c,d]}(N,\K,\eta)$.
When $c = -\infty$ and $d = +\infty$ then it is an invariant up to Liouville deformation.

To every Liouville domain $N$ we can form its completion $\widehat{N}$ by attaching
a cylindrical end $[1,\infty) \times \partial N$ to $\partial N$ and extending $\theta_N$
by $r_N d\alpha_N$ where $\alpha_N = \theta_N|_{\partial N}$ and $r_N$ is the coordinate
parameterizing $[1,\infty)$ called the cylindrical coordinate.
A Hamiltonian $H : S^1 \times \widehat{N} \rightarrow \R$
is said to be {\it admissible} if $H(t,x) = \lambda r_N(x)$ near infinity where $\lambda$
is a constant called the {\it slope} of $H$.
We sometimes view $H$ as a family of Hamiltonians $H_t : \widehat{N} \rightarrow \R$ where $t \in S^1$.
We have an $S^1$ family of vector fields $X_{H_t}$ and it
has an associated flow $\Phi^t_{X_{H_t}}$
(a family of symplectomorphisms parameterized by $t \in \R$ satisfying
$\frac{\partial}{\partial t} \Phi^t_{X_{H_t}} = X_{H_t}$
where we identify $S^1 = \R / \Z$).
A $1$-periodic orbit $o : S^1 \rightarrow \widehat{N}$
is a map which satisfies $o(t) = \Phi^t_{X_{H_t}}(x)$ for some $x \in \widehat{N}$.
We say that $o$ is non-degenerate if
$D\Phi^1_{X_{H_t}} : T_x\widehat{N} \rightarrow T_x\widehat{N}$
has no eigenvalue equal to $1$.
First of all, we can perturb $H$ slightly so that its slope $\lambda$
is not in the period spectrum ${\mathcal P}_N$.
This means that all of its $1$ periodic orbits sit inside some compact subset of $\widehat{N}$.
We then perturb $H$ again by a $C^\infty$ small amount so that
all of its $1$-periodic orbits are non-degenerate
and so that it still remains admissible (see \cite[Theorem 9.1]{SalamonZehnder:Morse}).
Because we have a trivialization $\tau$ of the canonical bundle
of $N$,
this gives us a canonical
trivialization of the symplectic bundle $TN$
restricted to an orbit $o$ (up to homotopy). Using this trivialization, we can define an index of $o$ called the
{\bf Robbin-Salamon} index (this is equal to the Conley-Zehnder
index taken with negative sign).
We will write $i(o)$ for the index of this orbit $o$.
For a 1-periodic orbit $o$ we define the {\bf action} $A_H(o)$ as:
\[A_H(o) := -\int_0^1 H(t,o(t))dt -\int_o \theta_N.\]

Choose a coefficient field $\K$ and an $S^1$ family of
almost complex structures $J_t$ compatible with
the symplectic form.
We assume that
$J_t$ is convex with respect to this cylindrical end outside some large compact set
(i.e. $\theta \circ J_t = dr$). We also say that $J_t$ is {\bf admissible} if such a condition holds.
Let
\[CF_k^d(H,J,\eta) := \bigoplus_{o} \K \langle o \rangle\]
where we sum over $1$-periodic orbits $o$ of $H$ satisfying
$A_H(o) \leq d$ whose Robbin-Salamon index is $k$.
We write 
\[CF_k^{(c,d]}(H,J,\eta) := CF_k^d(H,J,\eta) / CF_k^c(H,J,\eta).\]
As a vector space, $CF_k^{(c,d]}(H,J,\eta)$ does not depend on $J$ or $b$,
but the differential will.
We need to define a differential for the chain
complex $CF_k^d(H,J,\eta)$ such that the inclusion
maps $CF_k^c(H,J,\eta) \hookrightarrow CF_k^d(H,J,\eta)$ for $c<d$ are chain maps.
This makes $CF_k^{(c,d]}(H,J,\eta)$ into a chain complex as well.

We will now describe the differential \[\partial : CF_k^d(H,J,\eta) \rightarrow CF_{k-1}^d(H,J,\eta).\]
We consider curves
$u : \R \times S^1 \longrightarrow \widehat{N}$ satisfying
the perturbed Cauchy-Riemann equations:
\[ \partial_s u + J_t \partial_t u = \nabla^{g_t} H\]
where $\nabla^{g_t}$ is the gradient associated to the $S^1$ family of metrics
$g_t := \omega(\cdot,J_t(\cdot))$.
For two $1$ periodic orbits $o_{-},o_{+}$ let
$\overline{U}(o_{-},o_{+})$ denote the set of all curves $u$ satisfying
the perturbed Cauchy-Riemann equations such that $u(s,\cdot)$ converges
to $o_{\pm}$ as $s \rightarrow \pm \infty$. This has a natural
$\R$ action given by translation in the $s$ coordinate.
Let $U(o_{-},o_{+})$ be equal to $\overline{U}(o_{-},o_{+}) / \R$.
For a $C^{\infty}$ generic admissible complex structure
we have that $U(o_{-},o_{+})$ is an $i(o_{-}) - i(o_{+}) -1$  dimensional
manifold (see \cite{FHS:transversalitysymplectic}).
There is a maximum principle to ensure
that all elements of $U(o_{-},o_{+})$ stay inside
a compact set $K$ (see \cite[Lemma 1.5]{Oancea:survey} or
\cite[Lemma 7.2]{SeidelAbouzaid:viterbo}).
Hence we can use a compactness theorem (see for instance
\cite{BEHWZ:compactnessfieldtheory}) to ensure that
if $i(o_{-}) - i(o_{+}) = 1$, then
$U(o_{-},o_{+})$ is a compact zero dimensional manifold.
The class $b \in H^2(N,\Z / 2\Z)$ enables us to orient this manifold
(see \cite[Section 3.1]{Abouzaid:contangentgenerate})).
Let $\# U(x_{-},x_{+})$ denote the number of
positively oriented points of $U(x_{-},x_{+})$
minus the number of negatively oriented points. Then we have a differential:
\[\partial : CF_k^d(H,J,\eta) \longrightarrow CF_{k-1}^d(H,J,\eta) \]
\[\partial \langle o_{-} \rangle := \displaystyle \sum_{i(o_{-}) - i(o_{+})=1 } \# U(o_{-},o_{+}) \langle o_{+} \rangle\]
By analyzing the structure of 1-dimensional moduli spaces, one shows
$\partial^2=0$ and defines
$HF_*(H,\eta)$ as the homology of the above chain complex.
The homology group $HF_*^d(H,\eta)$ depends on $H$ and $\eta$ but is independent of $J$
up to canonical isomorphism.
We define $HF_*^{(c,d]}(H,\eta)$ as
the homology of the chain complex \[CF_*^d(H,J,\eta) / CF_*^c(H,J,\eta).\]

If we have two non-degenerate admissible
Hamiltonians $H_1 < H_2$, then there is a natural
map:
\[HF_*^{(c,d]}(H_1,\eta) \longrightarrow HF_*^{(c,d]}(H_2,\eta)\]
This map is called a {\bf continuation map}.
This map is defined from a map $C$ on the chain level as follows:
\[C : CF_k^d(H_1,J_1,\eta) \longrightarrow CF_k^d(H_2,J_2,\eta) \]
\[\partial \langle o_{-} \rangle := \displaystyle
\sum_{i(o_{-}) = i(o_{-}) } \# P(o_{-},o_{+}) \langle o_{+} \rangle\]
where $P(o_{-},o_{+})$ is a compact oriented zero dimensional manifold
of solutions of the following equations:
Let $K_s$, $s \in \R$ be a smooth non-decreasing family of admissible Hamiltonians 
equal to $H_1$ for $s \ll 0$ and $H_2$ for $s \gg 0$
and $J_{s,t}$ a smooth family of admissible almost complex structures joining $J_{1,t}$ and $J_{2,t}$.
The set  $P(o_{-},o_{+})$ is the
set of solutions to the parameterized Floer equations
\[ \partial_s u + J_{s,t}\partial_t u = \nabla^{g_t} K_{s,t}\]
such that $u(s,\cdot)$ converges
to $o_{\pm}$ as $s \rightarrow \pm \infty$.
For a $C^{\infty}$ generic family $(K_s,J_s)$ this is a compact zero dimensional
manifold.
Again the class $b \in H^2(N,\Z / 2\Z)$ enables us to orient this manifold.
If we have another such non-decreasing family admissible Hamiltonians joining $H_1$ and $H_2$
and another smooth family of admissible almost complex structures joining $J_1$ and $J_2$,
then the continuation map induced by this second family is chain homotopic to
the map induced by $(K_s,J_s)$.
The composition of two continuation maps is a continuation map.
If we take the direct limit of all these maps with respect
to admissible Hamiltonians $H$ ordered by $<$
such that $H|_N < 0$, then
we get our symplectic homology groups $SH_*^{(c,d]}(N,\eta)$.
We will write $SH_*^{\#}(N,\eta)$ for $SH_*^{(0,\infty)}(N,\eta)$.

Also we will write $SH_*$ instead of $SH_*^{(-\infty,\infty)}$.
If we wish to stress which coefficient field we are using,
we will write $SH_*^{\#}(M,\eta,\K)$ if the field is $\K$
for instance.
We will write $SH_*^{\leq d}$ instead of $SH_*^{(-\infty,d]}$.
We will suppress the term $\eta$ from the notation when the context is clear.
Also from now on whenever we have a Liouville domain or symplectic manifold then
we will assume that we have chosen such a pair $\eta = (\tau,b)$.

\subsection{Growth rates} \label{section:growthrates}

In order to define growth rates, we will need some linear algebra first.
Let $(V_x)_{x \in [1,\infty)}$ be a family of vector spaces indexed by $[1,\infty)$.
For each $x_1 \leq x_2$ we will assume that there is a homomorphism
$\phi_{x_1,x_2}$ from $V_{x_1}$ to $V_{x_2}$ with the property that
for all $x_1 \leq x_2 \leq x_3$, $\phi_{x_2,x_3} \circ \phi_{x_1,x_2} = \phi_{x_1,x_3}$
and $\phi_{x_1,x_1} = \text{id}$.
We call such a family of vector spaces a {\it filtered directed system}.
Because these vector spaces form a directed system, we can take the direct
limit $V := \varinjlim_x V_x$.
From now on we will assume that $V_x$ is finite dimensional.
For each $x \in [1,\infty)$ there is a natural map:
\[q_x : V_x \rightarrow \varinjlim_x V_x. \]
Let $a : [1,\infty) \rightarrow [0,\infty)$ be a function such that
$a(x)$ is the rank of the image of the above map $q_x$.
We define the growth rate as:
\[\Gamma( (V_x) ) : =\varlimsup_x \frac{\log{a(x)}}{\log{x}} \in \{-\infty\} \cup [0,\infty].\]
If $a(x)$ is $0$ then we just define $\log{a(x)}$ as $-\infty$.
If $a(x)$ was some polynomial of degree $n$ with positive leading coefficient, then the growth rate
would be equal to $n$.
If $a(x)$ was an exponential function with positive exponent, then the growth rate is $\infty$.

In the previous section we defined for a Liouville domain $N$
(whose boundary had discrete period spectrum), $SH_*^{\leq \lambda}(N)$.
For $\lambda_1 \leq \lambda_2$, there is a natural map
$SH_*^{\leq \lambda_1}(N) \rightarrow SH_*^{\leq \lambda_2}(N)$
given by inclusion of the respective chain complexes.
This is a filtered directed system $(SH_*^{\leq \lambda}(N))$
whose direct limit is $SH_*(N)$.
\begin{defn} \label{defn:growthrate}
We define the growth rate $\Gamma(N,\eta)$ as:
\[\Gamma(N,\eta) := \Gamma(SH_*^{\leq \lambda}(N,\eta))\]
\end{defn}
We also have the following theorem (\cite[Theorem 2.4]{McLean:affinegrowth}):
\begin{theorem} \label{thm:growthrateinvariance}
Let $N_1,N_2$ be two Liouville domains such that $\widehat{N_1}$
is symplectomorphic to $\widehat{N_2}$ where the symplectomorphism
pulls back $b_2 \in H^2(N_2,\Z / 2\Z)$ to $b_1 \in H^2(N_1,\Z / 2\Z)$
and $\tau_2$ to $\tau_1$ where $\tau_2$ and $\tau_1$ are trivializations of the canonical bundle.
Then
$\Gamma(N_1,(\tau_1,b_1)) = \Gamma(N_2,(\tau_2,b_2))$.
\end{theorem}

Hence we will just write $\Gamma(\widehat{N},d\theta_N,(\tau,b))$ for the growth rate of $(N,\theta_N)$.
We will sometimes just write $\Gamma(\widehat{N})$ if the context makes it clear
that $d\theta_N$ is our symplectic form and $(\tau,b)$ is our associated
trivialization and homology class.

\subsection{Local Floer homology}

In this section we mildly generalize the notion of local Floer homology as defined
in \cite{GinzburgGurel:localfloer}.
All the lemmas in this section and properties proven are almost exactly the same
as ones proven in \cite{GinzburgGurel:localfloer}.
Usually local Floer homology is defined for isolated $1$-periodic orbits
(see \cite{GinzburgGurel:localfloer}).
In our case we will define it for isolated families of $1$-periodic orbits which are all contained
inside some compact set and such that they have the same action.
Let $(Q,\omega_Q)$ be a symplectic manifold and $H : S^1 \times Q \rightarrow \R$ a Hamiltonian.
Let ${\mathcal F} \subset Q$ be a set of fixed points of $\phi^1_H$ inside $Q$.
We say that they are {\it isolated} if there is some open neighbourhood
${\mathcal N}_{\mathcal F}$ of ${\mathcal F}$ whose closure is compact
such that any fixed point of $\phi^1_H$ inside the closure of ${\mathcal N}_{\mathcal F}$
is contained inside ${\mathcal F}$.
We call ${\mathcal N}_{\mathcal F}$ an {\it isolating neighbourhood}.
Note that the orbits starting inside ${\mathcal F}$ can exit this neighbourhood,
we just require that they start at ${\mathcal F}$.

Let ${\mathcal F}$ be a set of fixed points of $\phi^1_H$ which is isolated and
such that the associated orbits have the same action.
We will now define a Floer homology group $HF_*(H_t,{\mathcal F})$
called Local Floer homology.
We need a lemma first.

\begin{lemma} \label{lemma:compactnessresult}
Let $G^n_t$ be a sequence of time dependent Hamiltonians which $C^\infty$
converge to $H_t$.
Let $J^n_t$ be a sequence of compatible almost complex structures $C^\infty$
converging to a compatible almost complex structure $J_t$.
Let $V$ be any open subset containing ${\mathcal F}$ whose closure is compact.
Let $U' \subset {\mathcal N}_{\mathcal F}$ be an open subset such that the flow
$\phi^t_{H_t}(U')$ is well defined for all $0 \leq t \leq 1$
(i.e. none of these points flow off to infinity).
Then for large enough $n$,
\begin{enumerate}
\item \label{item:orbitcompactness}

All $1$-periodic orbits $o(t)$ of $G^n_t$ starting inside ${\mathcal N}_{\mathcal F}$
must satisfy $o(t) \subset \phi^t_{H_t}(U')$ for $0 \leq t \leq 1$.
\item \label{item:floercompactness}

If $u : \R \times S^1 \rightarrow V$ is a Floer trajectory with respect to $(G^n_t,J^n_t)$
connecting orbits of $G^n_t$ starting inside ${\mathcal N}_{\mathcal F}$
then $u(s,t) \subset \phi^t_{H_t}(U')$ for all $0 \leq t \leq 1$ and $s \in \R$.
\end{enumerate}
\end{lemma}
\proof of Lemma \ref{lemma:compactnessresult}.
We identify $S^1 = \R / \Z$.
Suppose for a contradiction there is a subsequence $n_i$, a sequence of orbits $o_i$ of
$G^{n_i}_t$ starting inside ${\mathcal N}_{\mathcal F}$ and a sequence of points $t_i \in S^1$ so that
$o_i(t_i) \notin \phi^{t_i}_{H_{t_i}}(U')$.
By passing to a subsequence we can assume that $t_i$ converges to some point $t \in S^1$
and that the starting point $o_i(0)$
converges to some point $p$ in the closure of ${\mathcal N}_{\mathcal F}$.
Hence $o_i$ converges to some orbit $o$ of $H_t$ in the $C^0$ sense
and so $o_i(t_i)$ converges to $o(t)$.
Because $o$ is an orbit starting at $p$ which is contained inside the closure of
${\mathcal N}_{\mathcal F}$, we have that $p \in {\mathcal F}$.
Hence $o(t) \in \phi^t_{H_t}(U')$ for all $t \in [0,1]$ which is impossible because $o_i(t_i) \notin \phi^{t_i}_{H_t}(U')$ for all $i$.
Hence for large enough $n$ we have shown that
all $1$-periodic orbits $o(t)$ of $G^n_t$ starting inside ${\mathcal N}_{\mathcal F}$
must satisfy $o(t) \subset \phi^t_{H_t}(U')$ when $0 \leq t \leq 1$.

Suppose for a contradiction there is a sequence of Floer trajectories
\[ u_i : \R \times S^1 \rightarrow V \]
 with respect to
$(G^{n_i}_t,J^{n_i}_t)$
connecting orbits of $G^{n_i}_t$ starting inside $U'$ and a sequence of points $(s_i,t_i)$
such that $u_i(s_i,t_i) \notin \phi^{t_i}_{H_t}(U')$.
We can also assume that $u(s_i,0)$ is contained inside ${\mathcal N}_{\mathcal F}$.
The point is that if $u_i(s_i,0)$ was not contained inside this open set
for infinitely many $i$ then we know that $u_i(s,0) \in {\mathcal N}_{\mathcal F}$
for $s$ large enough (because $u_i$ converges to orbits starting inside ${\mathcal N}_{\mathcal F}$),
so by the continuity of $u_i$ we could find another sequence of points
$(s'_i,t'_i)$ with $t'_i = 0$ and
$u_i(s'_i,t'_i) \in {\mathcal N}_{\mathcal F} \setminus \phi^{t_i}_{H_t}(U')$.

We would like to use a compactness argument
(such as \cite{BEHWZ:compactnessfieldtheory})
to say that these Floer trajectories must converge to some Floer trajectory of energy $0$
but the problem is that the Hamiltonian $H$ could be degenerate.
So instead we will do the following:
First of all after shifting in the $s$ coordinate, we may as well assume that $s_i = 0$
for all $i$ and after passing to a subsequence we can assume that $t_i$ converges to $t'$
and $u_i(0,0)$ converges to some point $p$.
We have that $p$ is contained in the closure
of ${\mathcal N}_{\mathcal F}$.
We can view the maps $u_i$ as a sequence of holomorphic sections of a $C^\infty$
converging family of Hamiltonian fibrations
(see \cite{McduffSalamon:jholomorphic}) whose fiber is $Q$.
These fibrations converge to ${\mathcal H}$
which is the Hamiltonian fibration over $\R \times S^1$ associated to $H_t$.
Hence by using a compactness result such as
\cite{fish:compactness}  we have that for every compact subsurface $S$ of $\R \times S^1$,
 $u_i|_S$ converges in the Gromov sense
to some nodal curve $u : S' \rightarrow {\mathcal H}$.
Some of the components of this nodal curve could be holomorphic maps into the fibers of ${\mathcal H}$
(bubbles) and others are multisections.
Also these bubbles have energy $0$ and hence must be points inside the fibers.
There is at most one multisection $\tilde{u}_S$ and this must be a section
because our nodal curve intersects each fiber with multiplicity $1$.
By viewing $\R \times S^1$ as a union of compact surfaces
$\{-i \leq s \leq i\}$ we get after passing to a subsequence
and using the above compactness argument a section $\tilde{u}$
of ${\mathcal H}$.
We view this section $\tilde{u}$ as a map $u : \R \times S^1 \rightarrow Q$
satisfying the Floer equations.
This section has the property that $u(0,0) = p$.
The map $u$ has zero energy and hence $\frac{\partial u}{\partial s} = 0$
and $\frac{\partial u}{\partial t} = X_{H_t}$. This means that $u(0,t)$
is an orbit of $H$ starting at $p$ but this is impossible because
$u(0,t') \notin \phi^1_{H_{t'}}(U')$.
Hence for large enough $n$, $u(s,t) \subset \phi^t_{H_t}(U')$ for all $s \in \R$ and $0 \leq t \leq 1$.
\qed

We will now define $HF_*(H_t,{\mathcal F})$.
We choose some relatively compact open set $W$ containing all the orbits starting at ${\mathcal F}$
and an isolating neighbourhood ${\mathcal N}_{\mathcal F}$ for ${\mathcal F}$ whose closure is a subset of $W$.
We perturb $H_t$ very slightly to $H'_t$ so that all of its orbits are non-degenerate.
Choose a regular $S^1$ family of compatible almost complex structures $J_t$.
By the above Lemma we can ensure that all orbits starting inside ${\mathcal N}_{\mathcal F}$
are contained inside $W$ and the Floer trajectories
with respect to $(H'_t,J_t)$
connecting them inside $W$ are also contained inside $W$ (this is because we can choose $U'$
so that $\cup_{t \in [0,1]} \phi^t_{H_t}(U') \subset W$).
Also if a Floer trajectory connecting these orbits breaks then each component
must converge to an orbit starting at ${\mathcal N}_{\mathcal F}$ by this Lemma
because we can choose $U'$ so that its closure is contained in ${\mathcal N}_{\mathcal F}$.
Hence we have a well defined differential on the Floer chain complex generated by these orbits.
We define $HF_*(H_t,{\mathcal F})$ to be the homology of the Floer complex defined
using these orbits and Floer trajectories.

\begin{lemma} \label{lemma:independentofchoices}
We have that $HF_*(H_t,{\mathcal F})$ does not depend on
the choice of $(H'_t,J_t)$  or isolating neighbourhood
as long as the perturbation $H'_t$ is sufficiently small.
This also means that if we have some symplectomorphism $\phi$
from some neighbourhood ${\mathcal N}_{\mathcal F}$
to another isolating neighbourhood ${\mathcal N}'_{\mathcal F'}$
coming from some Hamiltonian $K_t$ so that $\phi^* K_t = H_t$ then
$HF_*(H_t,{\mathcal F}) = HF_*(K_t,\phi({\mathcal F}))$.
\end{lemma}
Note that these groups do depend on the choice of trivialization of the canonical bundle
and of the choice of class $b \in H^2(W,\Z / 2\Z)$.
But in the cases that we will use, the neighbourhood $W$ is homotopic to a $1$ complex
so $b$ must be zero.
\proof of Lemma \ref{lemma:independentofchoices}.
Let $(H'',J''_t)$ be another pair and let ${\mathcal N}'_{\mathcal F}$ another neighbourhood.
Let $(K_s,Y_s)$ be a smooth family of pairs parameterized by $s \in \R$
such that $(K_s,Y_s) = (H',J_t)$ for $s$ very negative and $(K_s,Y_s) = (H'',J''_t)$
for $s$ very positive.
Choose any relatively compact open set $W$ containing the orbits starting at ${\mathcal F}$.
If the perturbations $H''$ and $H'$ are small enough,
then the Floer trajectories $u(s,t)$ for $(H',J_t)$ and $(H'',J''_t)$
are contained inside an arbitrarily small open subset containing the orbits
by Lemma \ref{lemma:compactnessresult}.
Also if $K_s$ is sufficiently $C^\infty$ close to $H$ for all $s$, 
we get (by using the same proof as in Lemma \ref{lemma:compactnessresult})
that the continuation map Floer trajectories $u(s,t)$ for $(K_s,Y_s)$ inside $W$
are contained inside an arbitrarily small open subset containing the orbits.
This means our chain complexes are independent of the choice of isolating neighbourhood
and we have well defined continuation maps between them.
Hence we can use continuation arguments (for instance from
\cite[Section 6]{SalamonZehnder:Morse}) to prove invariance of choices
of $(H',J_t)$ and neighbourhood ${\mathcal N}_{\mathcal F}$.
\qed

Let $H^s_t$ be a smooth family of time dependent Hamiltonians parameterized by $s \in [0,1]$.
Let ${\mathcal F} \subset Q$ be an isolated set of fixed points of $H^s_t$ for every $s$.
If there is some isolating neighbourhood ${\mathcal N}_{\mathcal F}$ (independent of $s$) of these fixed points
for each $s \in [0,1]$ then we say that $(H^s_t,{\mathcal F})$ is an {\it isolated deformation}.

\begin{lemma}  \label{lemma:deformationinvariance}
Suppose we have an isolated deformation $(H^s_t,{\mathcal F})$ then
\[HF_*(H^0_t,{\mathcal F}) = HF_*(H^1_t,{\mathcal F}).\]
\end{lemma}
\proof of Lemma \ref{lemma:deformationinvariance} (Sketch).
First of all we add a smooth family of constants to $H^s_t$ so that
all the orbits starting at ${\mathcal F}$ have the same action for each $s$.
By using similar compactness ideas from Lemma \ref{lemma:compactnessresult},
we get a well defined continuation map from $HF_*(H^0_t,{\mathcal F})$
to $HF_*(H^1_t,{\mathcal F})$.
This is an isomorphism on homology as it has an inverse and continuation
maps are functorial (\cite[Section 6]{SalamonZehnder:Morse}).
\qed

Let $V'$ be a symplectic manifold and 
let $V$ be a codimension $0$ connected symplectic submanifold
and let $G_t$ be a Hamiltonian such that
$\phi^1_{G_t}$ is the identity map on $V$.
We suppose that $V'$ has a trivilization of its canonical bundle. 
Let $p \in V$ and consider the loop:
\[l(s) = \phi^s_{G_t}(p).\]
The choice of trivialization of the canonical bundle gives us a canonical
trivialization up to homotopy of the symplectic bundle $l^* (TV')$
which we view as a map $\tau$ from $l^* (TV')$ to $\C^n \cong T_{l(0)} V'$.
Hence we have a loop of linear symplectic automorphisms of $T_{l(0)} V' \cong \C^n$
given by $\tau \circ D\phi^s_{G_t}$.
This has an associated Maslov index $\mu$.
We say that $\mu$ is the Maslov index of the Hamiltonian loop generated by $G_t$.
If we have a Hamiltonian $H_t$ defined on $V'$ then we can form a new Hamiltonian
$(H \# G)_t$ by first modifying $H_t$ and $G_t$ so that near $t = 0,1$
these Hamiltonians are zero (by multiplying them by an appropriate function $\rho(t)$).
We then define $(H \# G)_t$ to be $H_{2t}$ for $0 \leq t \leq \frac{1}{2}$
and $G_{2t-1}$ for $\frac{1}{2} \leq t \leq 1$.

\begin{lemma} \label{lemma:effectofcircleaction}
Let ${\mathcal F}$ be an isolated set of fixed points of $H_t$ of the same action.
Suppose that we have a family of Hamiltonians $G_t$ such that the time $1$ flow is the identity map on a connected neighbourhood $V$ of our isolated orbits starting at ${\mathcal F}$. 
Then $HF_*(H_t,{\mathcal F}) = HF_{*-2\mu}( (H \# G)_t,{\mathcal F}')$ where $\mu$ is the Maslov index of our action $G_t$.
\end{lemma}
We will omit the proof of this Lemma as the key ideas are contained in
\cite[Section 2.3]{Ginzburg:Conley}.
This Lemma basically says that local Floer homology only depends on the time $1$-flow
of our Hamiltonian symplectomorphism locally around ${\mathcal F}$
up to some shift in index.
Note that we really need the compactness result \ref{lemma:compactnessresult}
to ensure that the orbits $o(t)$ of $(H \# G)_t$ and Floer trajectories $u(s,t)$
stay near $\phi^t_{(H \# G)_t}(V)$.

\begin{lemma} \label{lemma:disjointunionisolated}
Suppose that ${\mathcal F}$ and ${\mathcal F}'$
are two isolated fixed point sets of $H_t$ whose union is also an isolated fixed point set.
Then:
\[ HF_*(H_t,{\mathcal F} \cup {\mathcal F}') \cong
HF_*(H_t,{\mathcal F}) \oplus HF_*(H_t,{\mathcal F}').\]
\end{lemma}
\proof of Lemma \ref{lemma:disjointunionisolated}.
We can choose a small isolating neighbourhood of
${\mathcal F} \cup {\mathcal F}'$ 
which is the disjoint union of two isolating
neighbourhoods ${\mathcal N}_{\mathcal F}$
and ${\mathcal N}_{\mathcal F'}$.
For a small enough perturbation of $H_t$
we have that all the orbits $o(t)$ and Floer trajectories $u(s,t)$
that are used to define
$HF_*(H_t,{\mathcal F} \cup {\mathcal F}')$
satisfy $o(0), u(s,0) \in {\mathcal N}_{\mathcal F} \cup {\mathcal N}_{\mathcal F}'$ for all $s \in \R$.
Hence there are no Floer trajectories connecting
orbits starting inside ${\mathcal N}_{\mathcal F}$
with orbits starting inside ${\mathcal N}_{\mathcal F'}$.
Hence the chain complex defining
$HF_*(H_t,{\mathcal F} \cup {\mathcal F}')$
is the direct sum of the chain complexes defining
$HF_*(H_t,{\mathcal F})$ and $HF_*(H_t,{\mathcal F}')$.
This gives us our result.
\qed

Let $H_t$ be a time dependent Hamiltonian on $Q$.
The action spectrum of $H_t$ is the set of action values of all its $1$ periodic orbits.
\begin{lemma} \label{lemma:levelset}
Let $c \in \R$ and let ${\mathcal F}$ be the set of fixed points of $H_t$ of action $c$.
Suppose that the action spectrum  of $H_t$ is discrete in a neighbourhood of $c$
and ${\mathcal F}$ is compact.
Then ${\mathcal F}$ is an isolated family of orbits and for $\delta>0$ small enough,
\[HF_*(H_t,{\mathcal F}) = HF_*^{(c-\delta,c+\delta)}(H_t).\]
\end{lemma}
\proof of Lemma \ref{lemma:levelset}.
We can choose $\delta>0$ small enough so that the only orbits
of action in $(c-2\delta,c+2\delta)$ have action exactly $c$.
Any orbit starting at a point near ${\mathcal F}$
has action near $c$. But this means that this orbit has action $c$ and so
this orbit starts inside ${\mathcal F}$.
Hence ${\mathcal F}$ is isolated.
In order to define $HF_*^{(c-\delta,c+\delta)}(H_t)$,
we perturb $H_t$ slightly to a non-degenerate Hamiltonian $H'_t$
and then build our Floer complex using orbits only in the action window
$(c-\delta,c+\delta)$.
If we choose a small enough perturbation $H'_t$ of $H_t$ all of whose orbits are non-degenerate
then all the orbits of action $(c-\delta,c+\delta)$
are contained inside our isolating neighbourhood ${\mathcal N}_{\mathcal F}$.
Hence the chain complexes defining $HF_*^{(c-\delta,c+\delta)}(H_t)$
and $HF_*(H_t,{\mathcal F})$ are identical.
\qed

By a spectral sequence argument we get the following corollary of
Lemmas \ref{lemma:disjointunionisolated} and \ref{lemma:levelset}:
\begin{corollary} \label{corollary:spectralsequenceupperbound}
Let $H$ be a Hamiltonian with the property that $HF_*(H)$ is well defined
(in our case $H$ will be some admissible Hamiltonian on the completion of a Liouville domain).
Suppose also that the fixed points of $H$ form a disjoint union
of isolated families ${\mathcal F}_i$ $i = 1,\cdots,l$.
Then the rank of $HF_k(H)$ is bounded above by the rank
of $\bigoplus_{i=1}^l HF_k(H,{\mathcal F}_i)$.
\end{corollary}

Let $p$ be an isolated fixed point of the Hamiltonian symplectomorphism
induced by $H_t$.
Then from \cite{GinzburgGurel:localfloer}, there is an index $\Delta_{H_t}(p) \in \R$ satisfying:
\begin{enumerate}
\item \cite[Property MI1]{GinzburgGurel:localfloer}.

$\Delta_{kH_{kt}}(p) = k \Delta_{H}(p)$.

\item \cite[Property LF5]{GinzburgGurel:localfloer}.

Let $n$ be half the dimension of our symplectic manifold.
Then $HF_k(H_t,p)$ is zero if $k \notin [\Delta_{H_t}(p) - n, \Delta_{H_t}(p) + n]$.
\end{enumerate}

\begin{lemma} \label{lemma:localkunneth}
Suppose we have a neighbourhood of the orbits starting at ${\mathcal F}$
which is symplectomorphic to a product $U \times V$ and where the
Hamiltonian splits up as $H_U + H_V$ where $H_U$ is a
Hamiltonian on $U$ and $H_V$ is a Hamiltonian on $V$.
Then $HF_*(H,{\mathcal F}) =HF_*(H_U,{\mathcal F}_1) \otimes HF_*(H_V,{\mathcal F}_2) $ where ${\mathcal F}  = {\mathcal F}_1 \times {\mathcal F}_2$.
The choice of trivialization of the canonical bundle and our choice of class $b$ must also split up as a product.
\end{lemma}
The reason why this Lemma is true is because we can perturb
our Hamiltonian so that it is still a product and also the choice of almost complex structure
can also split up as a product. This ensures that the chain complex splits up as a tensor product.
A very similar statement is contained in 
\cite[Property (LF4) in Section 3.2]{GinzburgGurel:localfloer}.

\subsection{Reeb orbit homology theories}

Let $M$ be a Liouville domain.
We choose some trivialization $\tau$ of its canonical bundle.
This induces a trivialization $\tau_\partial$
of the canonical bundle associated to the contact distribution on $\partial M$.
This is because the symplectic complement of the contact distribution is a
symplectic bundle trivialized by the Reeb vector field and the Liouville vector field.
Let $\gamma$ be a (not necessarily simple) Reeb orbit of $\partial M$.
We view $\gamma$ as a map from $\R / l_\gamma \Z$
to $\partial M$ so that $\frac{d}{dt}(\gamma(t)) = R$
where $R$ is the Reeb vector field.
Here $l_\gamma$ is the length of the Reeb orbit.
Note that $\gamma(t+c)$ is a Reeb orbit for any constant $c \in \R$.
We assume that $\gamma$ is {\it isolated}.
This means that there is some neighbourhood $U$ of $\gamma$
such that there are no Reeb orbits intersecting $U \setminus \text{image}(\gamma)$.
We can define an invariant $CH_*(\gamma)$ as follows:
Because $\frac{d}{dt}(\gamma(t))$ is in the kernal
of $d\alpha_M$, we can find a fibration
$\pi_\gamma : {\mathcal N}_\gamma \twoheadrightarrow S^1$
where ${\mathcal N}_\gamma$ is a small neighbourhood of $\text{image}(\gamma)$ and
such that $d\alpha_M$ restricted to each fibre is a symplectic form.
By possibly shrinking ${\mathcal N}_\gamma$ and using a Moser theorem,
we can assume each fibre is symplectomorphic to a small
ball $B_\delta \subset \R^{2n-2}$ of radius $\delta >0$ and the structure
group is $U(n-1)$.
Because $\pi_\gamma$ is a fibration, we have a vertical tangent bundle
(i.e. the subbundle of the tangent bundle which is tangent to the fibers of $\pi_\gamma$).
Consider the vertical tangent bundle restricted to the zero section of $\pi_\gamma$.
This is homotopic through symplectic bundles to the contact distribution.
Hence we trivialize $\pi_\gamma$ so that the highest exterior power
of the vertical tangent bundle along the zero section coincides with our trivialization $\tau_\partial$.
This choice of trivialization is unique up to homotopy.
Hence $\pi_\gamma$ is a product fibration
$S^1 \times B_\delta$ and $d\alpha_M$ restricted to each
fibre is the standard symplectic form on $B_\delta$.
The line field spanned by the Reeb vector field
gives us a symplectic connection on this fibre
bundle because the Reeb vector field is in the kernal of $d\alpha_M$.
We define an $S^1$ family of vector fields on $B_\delta$ as follows:
For $t \in S^1$, we have a unique fiber $F_t$ of $\pi_\gamma$ which intersects
$\gamma(t)$. Our trivialization gives us a symplectomorphism $F_t \cong B_\delta$
and a natural projection $S^1 \times B_\delta \twoheadrightarrow B_\delta$.
By abuse of notation we write $\frac{d}{dt}$ for the vector field on $S^1$ given by
$(\pi_\gamma \circ \gamma)_* (\frac{d}{dt})$.
Let $\widetilde{\frac{d}{dt}}$ be the unique horizontal lift of $\frac{d}{dt}$.
We define $V_t$ to be the projection of $-\widetilde{\frac{d}{d t}}|_{F_t}$
to $B_\delta$.
Because $i(-\widetilde{\frac{d}{d t}})d\alpha_M = 0$
we have that $V_t$ is an $S^1$ family of symplectic vector fields.
These also preserve the origin.
Because $B_\delta$ is contractible and $V_t$ vanishes at the origin,
$V_t$ is generated by an $S^1$ family of Hamiltonians $H^\gamma_t$ which fix $0$.
Because the Reeb orbit is isolated, we have that $H^\gamma_t$
has an isolated fixed point at $0 \in B_\delta$.
Hence we define $CH_*(\gamma) := HF_*(H^\gamma_t,0)$.
This is independent of choice of fibration $\pi_\gamma$
because if we had two such fibrations then
we can join then via a smooth family of such fibrations.
Associated to this smooth family of fibrations we have a smooth family
of Hamiltonians $H^s_\gamma$ with isolated
fixed points at $0$ and hence by Lemma \ref{lemma:deformationinvariance}
they all have the same local Floer homology group.
We will call $CH_*(\gamma)$ the {\it Reeb orbit homology of $\gamma$}.

The problem with Reeb orbit homology is that there is not a very obvious link
between this homology group and symplectic homology.
So we now give another Floer homology group associated to this Reeb orbit
which has a slightly more direct relationship with symplectic homology.
Suppose that $\gamma$ has length $l_\gamma$.
Let $r_M$ be the cylindrical coordinate in $\widehat{M}$.
Choose a Hamiltonian $H_t$ on $\widehat{M}$ so that
there is some $x > 1$ with
$H_t = h(r_M)$ near $\{x\} \times \text{image}(\gamma) \subset [1,\infty) \times \partial M$
where $h'(x) = l_\gamma$ and $h''(l_\gamma) > 0$.
Then $H$ has an isolated $S^1$ family of fixed points
$\{x\} \times \text{image}(\gamma)$.
If $H_t$ has the above properties then we say that $H$
is {\it admissible with respect to} $\gamma$.
We define Reeb orbit Floer homology to be:
$HF_*(H_t,\{x\} \times \text{image}(\gamma))$.

Suppose I have another Hamiltonian $H'_t$ which is admissible with respect to $\gamma$
so that it has an isolated $S^1$ family of fixed points $\{x'\} \times \text{image}(\gamma)$.
Then there is a smooth family of Hamiltonians $H^s_t$ joining $H_t$ and $H'_t$
which are all admissible with respect to $\gamma$.
Hence $HF_*(H_t,\{x\} \times \text{image}(\gamma))$
is equal to $HF_*(H'_t,\{x'\} \times \text{image}(\gamma))$
by Lemma \ref{lemma:deformationinvariance}.
Hence this Floer homology group is independent of the choice of Hamiltonians which
are admissible with respect to $\gamma$.
We will call this group the
{\it symplectic homology of} $\gamma$
and we will write $SH_*(\gamma)$.

\section{Symplectic homology of iterates of a Reeb orbit}

Let $\gamma$ be a simple Reeb orbit and $\gamma^k$ its $k$-fold iterate.
The aim of this section is to prove:
\begin{theorem} \label{theorem:reeborbithomologybound}
There is some constant $C$ depending on our Reeb orbit $\gamma$
such that the rank of $SH_*(\gamma^k)$ is bounded above by $C$.
Also we can assign an index $\Delta(\gamma) \in \R$ for each isolated Reeb orbit $\gamma$
such that $SH_i(\gamma^k)$ is zero if $i \notin [k \Delta(\gamma) - n +1, k \Delta(\gamma) + n]$.
Here $n$ is half the dimension of our Liouville domain $M$. 
\end{theorem}

Let $B_\delta$ be an open ball in $\R^{2n}$ of radius $\delta>0$.
Let $Q$ be a symplectic manifold with a choice of diffeomorphism to
$B_\delta \times ([0,1] \times S^1)$ and with an exact symplectic form $d\theta_Q$.
Let $\pi_Q : B_\delta \times ([0,1] \times S^1) \rightarrow [0,1] \times S^1$
be the natural projection map.
Suppose that $Q$, $\pi_Q$ and $\theta_Q$ satisfy:
\begin{enumerate}
\item $d\theta_Q|_{B_\delta \times \bracket{(s,a)}}$ is the standard symplectic form on $B_\delta$
for all $(s,a) \in [0,1] \times S^1$.
Here we have identified $Q$ with $B_\delta \times ([0,1] \times S^1)$ and from now on we will do this.
\item
We require that the tangent spaces to the submanifold
\[\{0\} \times ([0,1] \times S^1) \subset B_\delta \times ([0,1] \times S^1)\]
are symplectically orthogonal to the fibers. 
\item
We let $\widetilde{\frac{\partial}{\partial a}}$ be the lift of
$\frac{\partial}{\partial a}$ up to the plane distribution ${\mathcal H}$
which is $d\theta_Q$ orthogonal to the fibers.
The vector field $\widetilde{\frac{\partial}{\partial a}}$ has no orbits
of any period contained inside ${\mathcal N}_Q \setminus \left( \{0\} \times ([0,1] \times S^1) \right)$
where ${\mathcal N}_Q$ is a small neighbourhood of $\{0\} \times ([0,1] \times S^1)$.
\item The symplectic form $d\theta_Q$ restricted to $\{0\} \times ([0,1] \times S^1)$
must be equal to
$\kappa_Q ds \wedge da$ for some constant $\kappa_Q >0$.
\end{enumerate}
We call such a fibration $\pi_Q$ a {\it partially trivial fibration}.
A {\it deformation of partially trivial fibrations} is a smooth family
of such fibrations where the map $\pi_Q$ is fixed along with the trivialization
$B_\delta \times ([0,1] \times S^1)$ but
the $1$-form $\theta_Q$ can smoothly vary and so can $\kappa_Q$.
The neighbourhood ${\mathcal N}_Q$ described above must be fixed throughout this deformation as well
although we are allowed to choose a smaller neighbourhood at the start of the deformation if we wish.

All such fibrations have a natural choice of trivialization of the canonical bundle
because the vertical bundle has a symplectic trivialization
induced by our choice of trivialization $B_\delta \times ([0,1] \times S^1)$
and the horizontal bundle has a symplectic trivialization induced by the coordinates $(s,a) \in [0,1] \times S^1$
where we view $S^1$ as the quotient $\R / \Z$ (so $da$ on $S^1$ has volume $1$).
Let $f : [0,1] \rightarrow \R$ be a function with $f',f'' >0 $
and $f'( \frac{1}{2} ) = \kappa_Q$.
Then ${\mathcal F}_Q := \{0\} \times (\{\frac{1}{2}\} \times S^1)$
is an isolated $S^1$ family of fixed points for $\pi_Q^*kf(s)$ for all $k \in \Z$.
We define $SH_*(\pi_Q,\theta_Q,k)$ to be
$HF_*(\pi_Q^*(kf(s)),{\mathcal F}_Q)$.
We project the vector field
$-\widetilde{\frac{\partial}{\partial a}}$
to a vector field $L$ tangent to the
fibers of $\pi_Q$ using the trivialization
$B_{\delta} \times ([0,1] \times S^1)$.
We view $L$ as a family of vector fields on
$B_{\delta}$ parameterized by $(s,a) \in [0,1] \times S^1$ .
These are Hamiltonian vector fields so
they are generated by a smooth family of
Hamiltonians $H^s_a$ which we will call
the {\it associated generating family of Hamiltonians
for }$Q$.
We define 
$CH_*(\pi_Q,\theta_Q,k)$ to be equal to
$HF_*(kH_{kt}^{\frac{1}{2}},0)$.
These groups are invariants of $Q$ up to deformation by Lemma 
\ref{lemma:deformationinvariance}.
We say that $\pi_Q$ is {\it trivial at infinity} if
$\theta_Q = \theta_B + C s da$ outside some compact subset of $Q = B_\delta \times ([0,1] \times S^1)$
for some constant $C$.
Here $\theta_B$ is a $1$-form on $B_\delta$ such that $d\theta_B$ is the standard
symplectic form on $B_\delta$.

\begin{lemma} \label{lemma:ensuringparalleltransportmapsarewelldefined}
Let $\pi_Q$ be a partially trivial fibration then it is deformation equivalent
to a partially trivial fibration which is trivial at infinity.
\end{lemma}
\proof of Lemma \ref{lemma:ensuringparalleltransportmapsarewelldefined}.
Let $\theta_B$ be a $1$-form on $B_\delta$ such that $d\theta_B$ is the standard
symplectic form on $B_\delta$.
We have that $\theta_Q = \theta_B + \beta + dR$ where $R$ is a function and
$\beta$ is a $1$-form which vanishes when restricted to the fibers of $\pi_Q$.
Let $\beta_t$ be a smooth family of $1$-forms such that $\beta_t = \beta$
near $\{0\} \times ([0,1] \times S^1)$ for all $t \in [0,1]$.
We also require that $\beta_0 = \beta$ and $\beta_1 = 0$
outside a small neighbourhood of $\{0\} \times [0,1] \times S^1$.
We have for a large enough constant $C>0$ that
\[ \theta^t_Q := \theta_B + \beta_t + d((1-t)R) + C t \pi_Q^* sda \]
is a deformation of partially trivial fibrations such that $\theta_Q = \theta^0_Q$.
Also because $\beta_1 = 0$ outside a small neighbourhood of $\{0\} \times ([0,1] \times S^1)$,
$(\pi_Q,\theta^1_Q)$ is trivial at infinity.
Hence we have a deformation of partially trivial fibrations starting at $\theta_Q$
and ending at one which is trivial at infinity.
\qed

Let $K_t : B_\delta \rightarrow \R$ be an $S^1$ family of compactly supported Hamiltonians.
We  suppose that $dK_t(0) = 0$ for all $t$, and that the constant $k$ periodic orbit at $0$ is isolated
for all $k \in \Z$.
We can construct a partially trivial fibration as follows:
We start with $B_{\delta} \times ([0,1] \times \R)$
with the product symplectic form $d\theta_{B_\delta} + ds \wedge da$.
This has a $\Z$ action where $1 \in \Z$ sends 
$(z,s,a) \in B_{\delta} \times ([0,1] \times \R)$
to $(\phi^{-1}_{K_t}(z),s,a+1)$.
We will define $Q_{K_t}$ to be the quotient
$B_{\delta} \times ([0,1] \times \R) / \Z$.
This has a trivialization 
\[T : B_{\delta} \times ([0,1] \times S^1)
\rightarrow Q_{K_t}\]
given by $T(z,s,a) = (\phi^{-a}_{K_a}(z),s,a)$.
Also $H^2(Q_{K_t})=0$ so the symplectic form has a primitive $\theta_{K_t}$.
We say that $(Q_{K_t},\theta_{K_t})$ is the {\it standard trivialization associated to} $K_t$.
Such fibrations are called {\it standard partially trivial fibrations}.

\begin{lemma} \label{lemma:standardformdeformation}
Every partially trivial fibration is deformation
equivalent to a standard partially trivial fibration.
\end{lemma}
\proof of Lemma \ref{lemma:standardformdeformation}.

First of all our partially trivial fibration
is deformation equivalent to some
partially trivial fibration $Q$ that is trivial at infinity.
The reason why we need a fibration trivial at infinity is that we have well defined
parallel transport maps (i.e. points don't get transported off to infinity).
Let $B_{\delta} \times ([0,1] \times S^1)$ be its choice of trivialization with respective coordinates
$(z,s,a)$. 
We have a family of smooth maps
$\psi_t : Q \rightarrow Q$ parameterized
 by $t \in [0,1]$ sending $(z,s,a)$
to $(z,\frac{1}{2}( 1 + (1-t)(2s-1)),a)$.
This is a smooth linear deformation retraction of $B_{\delta} \times ([0,1] \times S^1)$
onto $B_{\delta} \times (\{\frac{1}{2}\} \times S^1)$.
We define $\theta^t_q$ to be
$\psi_t^* \theta_Q$.
We have that $\psi_t^* \theta_Q$
is a symplectic form for $t<1$ but not for
$t=1$. But this problem can be fixed by
adding $Ct\pi_Q^* sda$ for some $C>0$.
So \[(Q,\theta^t_Q := \psi_t^* \theta_Q + Ct\pi_Q^* sda)\]
is a deformation of partially trivial fibrations.
The partially trivial fibration
$(Q,\theta^1_Q)$
has an associated family of Hamiltonians $H^s_a$ that are independent of $s$ so we will just write $H_a$.
These are all compactly supported.

On the trivialization
$B_{\delta} \times ([0,1] \times S^1)$
we have a smooth self diffeomorphism $\xi$ defined away from $B_\delta \times ([0,1] \times \{0\})$
given by sending
$(z,s,a)$ to $(\phi^{-a}_{H_a}(z),s,a)$ for $0 < a < 1$.
Let $\omega_{\text{std}}$
be the pullback
$\xi^*(d\theta_B + ds \wedge da)$ where $d\theta_B$ is the standard symplectic form on $B_\delta$.
This extends to a smooth form on $B_\delta \times ([0,1] \times S^1)$ which we define by abuse of notation as
$\omega_{\text{std}}$. There is a primitive $\theta_{\text{std}}$ such that $d\theta_{\text{std}} = \omega_{\text{std}}$.
This $1$ form gives $\pi_Q$
the structure of a partially trivial fibration in standard form along with the chosen trivialization 
$B_{\delta} \times ([0,1] \times S^1)$.
Also the horizontal lifts of $\frac{\partial}{\partial a}$
with respect to both 
$\theta^1_Q$ and
$\theta_{\text{std}}$
coincide.
Hence if $\rho : [0,1] \rightarrow \R$
is a smooth function which is zero at
$0$ and $1$ but positive elsewhere then for $\kappa > 0$ large enough we have that
\[(1-t) \theta^1_Q + t\theta_{\text{std}} + \kappa \rho(t) sda\]
is a deformation of partially trivial fibrations.
Hence $(Q,\theta_Q)$ is deformation equivalent to $(Q,\theta_{\text{std}})$ which is a standard partially trivial fibration.
\qed

\begin{lemma} \label{lemma:upperboundforpartialllytrivialfibration}
Let $(\pi_Q,\theta_Q)$ be a partially trivial fibration then the rank of $SH_l(\pi_Q,\theta_Q,k)$
is bounded above by the rank of \[CH_l(\pi_Q,\theta_Q,k) \oplus CH_{l-1}(\pi_Q,\theta_Q,k).\]
\end{lemma}
\proof of Lemma \ref{lemma:upperboundforpartialllytrivialfibration}.
We have that $SH_*(\pi_Q,\theta_Q,k)$ is a local Floer homology group
associated to an $S^1$ family of $1$-periodic orbits of some Hamiltonian.
In order to prove our lemma we will first deform our fibration $\pi_Q$
so that it is sufficiently nice.
We will then perturb our $S^1$ family of orbits so that they become two isolated orbits.
By analyzing these two isolated orbits we can relate them to $CH_*(\pi_Q,\theta_Q,k)$.

By Lemma \ref{lemma:standardformdeformation} we can assume that $(\pi_Q,\theta_Q)$
is a standard partially trivial fibration.
This has a universal cover which is a product
$B_{\delta} \times ([0,1] \times \R )$ with product symplectic form
$d\theta_{B_\delta}+ ds \wedge da$.
This also has an associated Hamiltonian $H_t$.
Let $f : [0,1] \rightarrow \R$ be a function with $f',f'' >0$ and $f'(\frac{1}{2})=1$.
Consider the function $kf(s)$ on $[0,1] \times S^1$ where $k \in \Z$.
This has an isolated $S^1$ family of fixed points $\{\frac{1}{2}\} \times S^1$
which are Morse Bott non-degenerate.
Let $\phi^t_s$ be the time $t$ flow of the Hamiltonian $s$.
The time $1$ flow is the identity map and this is a Hamiltonian $S^1$ action
of Maslov index $0$.
Let $\nu : S^1 \rightarrow \R$ be a Morse function with exactly one maximum and one minimum.
The Hamiltonian $K_t := kf(s) + \epsilon (\phi^{-kt}_s)^* \nu(a)$ is a small perturbation
of $kf(s)$ for $\epsilon > 0$ small enough.
The $S^1$ family of orbits $\{\frac{1}{2}\} \times S^1$ gets perturbed into
two orbits of index $0$ and $1$ corresponding to the maximum and minimum of $\nu$ respectively.
These points are located at $(\frac{1}{2},x_1)$ and $(\frac{1}{2},x_2)$ where $x_1$
and $x_2$ are the maximum and minimum points of $\nu$.
There is a small neighbourhood $U$ around the points $(\frac{1}{2},x_1)$ and $(\frac{1}{2},x_2)$ 
such that $\pi_Q^{-1}(U)$ is symplectomorphic to $U \times B_\delta$
with the standard product symplectic form and where $\pi_Q$ corresponds to the projection map to $U$.
Here the time $1$ flow of the Hamiltonian $\pi_Q^*(K_t)$ is equal to the time $1$
flow of the Hamiltonian $K'_t := \pi_Q^*(kf(s) - ks + \epsilon \nu(a)) + \pi_2^* kH_{kt}$
on $U \times B_\delta$ where $\pi_2$ is the natural projection map
$U \times B_\delta \twoheadrightarrow B_\delta$.
The reason for this is as follows:
If $p(t)$ is a path tangent to the vector field $X_1 := X_{kf(s) + \epsilon (\phi^{-kt}_s)^* \nu(a)}$
then $\phi^{-kt}_s(p(t))$ is tangent to the vector field $(\phi^{-kt}_s)_* X_1 + \frac{d}{dt}\phi^{-kt}_s$.
Hence by Lemma \ref{lemma:effectofcircleaction} we have
\[ HF_*\left(\pi_Q^* (K_t),\{0\} \times (\{\frac{1}{2}\} \times \{x_i\}) \right) =
HF_*\left( K'_t,\{0\} \times (\{\frac{1}{2}\} \times \{x_i\}) \right).\]
The point is that $(\pi_Q^* K_t \# (-K'_t))$ is a Hamiltonian $S^1$ action
on $\pi_Q^{-1}(U)$ isotopic through such actions to
the Hamiltonian $S^1$ action induced by $s$ which has Maslov index $0$.
The time $t$ flow of the Hamiltonian $K'_t$ fixes the points $(0,\frac{1}{2},x_1)$
and $(0,\frac{1}{2},x_2)$ on $U \times B_\delta$ for all $t$, so in particular
all the orbits starting at these points stay inside the product $U \times B_\delta$.
So by Lemma \ref{lemma:localkunneth} we get:
\[HF_*\left(K'_t,\{(0,\frac{1}{2},x_i)\} \right) =\]
\[HF_*\left(kf(s) - ks + \epsilon \nu(t),\{(\frac{1}{2},x_i)\}\right) \otimes HF_*\left( H_t,\{0\} \right).\]
Because the orbits of $kf(s) - ks + \epsilon \nu(a)$ are non-degenerate critical points
of index $0$ and $1$ we get that
\[HF_*\left(kf(s) -ks + \epsilon \nu(a),\{(\frac{1}{2},x_i)\}\right) = 
\left\{ \begin{array}{ll}
\K & \text{if } * = i \\
 0 & \text{otherwise}
\end{array}\right.\]
for $i = 0,1$.
Hence \[HF_*\left(K'_t,\{(0,\frac{1}{2},x_i)\} \right) =
HF_{*-i}\left( kH_{kt},\{0\} \right).\]
So by Lemma \ref{corollary:spectralsequenceupperbound} we have that the rank of
$HF_l(kf(s),\{0\} \times (\{\frac{1}{2}\} \times S^1) )$ is bounded above by the rank of
$\oplus_{i=0}^1 HF_{l-i}\left( kH_{kt},\{0\} \right)$.
Hence the rank of $SH_l(\pi_Q,\theta_Q,k)$
is bounded above by the rank of the group $CH_l(\pi_Q,\theta_Q,k) \oplus CH_{l-1}(\pi_Q,\theta_Q,k)$.
\qed

\proof of Theorem \ref{theorem:reeborbithomologybound}.
In this proof we will show that $SH_*(\gamma^k)$ is equal to $SH_*(\pi'_\gamma,\theta_M,k)$
for some partially trivial fibration $\pi'_\gamma$.
We will then use results from \cite{GinzburgGurel:localfloer} to put a bound
on $CH_*(\pi'_\gamma,\theta_M,k)$ and hence by Lemma
\ref{lemma:upperboundforpartialllytrivialfibration} we get our bounds on $SH_*(\gamma^k)$.

We will first assign an index to this Reeb orbit as follows:
There is a fibration map
$\pi_\gamma : {\mathcal N}_\gamma \twoheadrightarrow S^1$
where ${\mathcal N}_\gamma$ is a small neighbourhood of $\gamma$ and
such that $d\alpha_M$ restricted to each fibre is a symplectic form.
Here the fibers are symplectomorphic to the ball $B_\delta$ of radius $\delta$ for some $\delta > 0$.
This fibration also has a choice of trivialization which is compatible with the
trivialization of the canonical bundle on $\widehat{M}$.
If we look at $\{2\} \times R$ in the cylindrical end $[1,\infty) \times \partial M$
of $\widehat{M}$ then we have a fibration:
$\pi'_\gamma : [1,3] \times {\mathcal N}_\gamma \rightarrow [1,3] \times S^1$
where $\pi'_\gamma = (\text{id},\pi_\gamma)$.
This is a partially trivial fibration and we have that
$SH_*(\gamma^k) = SH_*(\pi'_\gamma,\theta_M,k)$
and $CH_*(\gamma,k) = CH_*(\pi'_\gamma,\theta_M,k)$.
By Lemma \ref{lemma:upperboundforpartialllytrivialfibration}
we have that the rank of $SH_l(\pi'_\gamma,\theta_M,k)$
is bounded above by the rank of $CH_l(\pi'_\gamma,\theta_M,k) \oplus CH_{l-1}(\pi_\gamma,\theta_M,k)$.
The fibration $\pi'_\gamma$ also has an associated family of Hamiltonians
$H^s_t$.
The Hamiltonian $kH^{\frac{1}{2}}_{kt}$ has an isolated fixed point at $0$
and so we can assign a mean index $\Delta_k := \Delta_{kH^{\frac{1}{2}}_{kt}}$.
We have by property (MI1) that $\Delta_k = k\Delta_1$.
We will define our index $\Delta(\gamma)$ to be $\Delta_1$.
Hence the rank of
$CH_l(\pi'_\gamma,\theta_M,k)$ is zero for $l \notin [k\Delta(\gamma) - (n-1), k\Delta(\gamma) + (n-1)]$
by property (LF5) stated above.
Hence $SH_*(\pi'_\gamma,\theta_M,k)$ is only supported in degrees
$[k\Delta(\gamma) - (n-1), k\Delta(\gamma) + n]$.
Also by \cite[Corollary 1.5]{GinzburgGurel:localfloer}
we get that the rank of $CH_*(\pi'_\gamma,\theta_M,k)$
is bounded above by some constant independent of $k$.
Hence the rank $SH_*(\gamma^k)$ is bounded above by some constant and independent of $k$.
\qed

\section{Proof of the main theorem}

Here is a statement of Theorem \ref{theorem:mainresult}:

{\it Suppose that $M$ is a Liouville domain such that $\partial M$
has only finitely many simple Reeb orbits, then:
\begin{enumerate}
\item There is a constant $C$ such that the rank of
$SH_k(M)$ is bounded above by $C$ for all $k \notin [1-n,n]$
where $n$ is half the dimension of $M$.
\item $\Gamma(M) \leq 1$.
\end{enumerate}
}

\proof of Theorem \ref{theorem:mainresult}.

Let $r$ be the radial coordinate on the cylindrical end $[1,\infty) \times \partial M$.
Let $h : [1,\infty) \rightarrow \R$ be a function which is $0$ near $1$ with
$h',h'' \geq 0$.
We also assume that $h'(x) = 1$ for $x \geq 2$ and that $h'' > 0$ in the region
where $0 < h' < 1$.
Let $\lambda \notin {\mathcal P}$ where ${\mathcal P}$ is the period spectrum of $\partial M$.
The Hamiltonian $\lambda h(r)$ on $\widehat{M}$ has the following isolated families of fixed points:
One family is the region $h^{-1}(0)$.
Also for each Reeb orbit $\gamma$ of length $l \leq \lambda$ there is a family
$\{\mathcal F\}_{\gamma,\lambda}$ equal to
$\{h_\gamma\} \times \gamma \subset [1,\infty) \times \partial M$
where $h_\gamma$ is the unique value that satisfies $\lambda h'(h_\gamma) = l$.
Because $h(r)$ is $C^2$ small in the region $h^{-1}(0)$,
we get $HF_*(h(r),h^{-1}(0)) = H^{n-*}(M)$.
Also by the definition of $SH_*(\gamma)$ we have that
$HF_*(h(r),{\mathcal F}_{\gamma,\lambda}) = SH_*(\gamma)$.
Hence by Corollary
\ref{corollary:spectralsequenceupperbound}, we get that the rank of
$HF_j(\lambda h(r))$ is bounded above by the rank of:
\[ H^{n-j}(M) \oplus \left( \oplus_{\gamma} SH_j(\gamma)\right)\]
where the direct sum $\oplus_\gamma$ is over all Reeb orbits of length $\leq \lambda$.

We have that symplectic homology is the direct limit as $\lambda$ tends to infinity
of $HF_*(\lambda h(r))$.
So the rank of $SH_j( M )$ is bounded above by the rank of
\[ H^{n-j}(M) \oplus \left(\oplus_{\gamma} SH_j(\gamma)\right)\]
where the sum $\oplus_{\gamma}$ is now over all Reeb orbits $\gamma$.
We have that $\partial M$ has only finitely many simple Reeb orbits $\gamma_1,\cdots,\gamma_m$.
We will write $\gamma^k_i$ for the $k$'th iterate.
By Theorem \ref{theorem:reeborbithomologybound},
there is a constant $C$ so that the rank of $SH_*(\gamma^k_i)$ is bounded above by
$C$ for all $k,i$.
Also we can assign an index $\Delta_i \in \R$ for each orbit $\gamma_i$
so that $SH_*(\gamma^k_i)$ is supported in degrees $[k\Delta_i - n+1,k\Delta_i + n]$.
This means that if $\Delta_i \neq 0$ then the rank of $\oplus_k SH_*(\gamma^k_i)$
is bounded in each degree.
If $\Delta_i = 0$ then $\oplus_k SH_*(\gamma^k_i)$ is supported only in degrees $[1-n,n]$.
Putting all of this together we get that the rank of $SH_k(M)$ is bounded above by some constant
independent of $k$ for all $k \notin [1-n,n]$.

We now need show that the growth rate is at most $1$.
Because the rank of $SH_*(\gamma^k_i)$ is bounded above by a constant,
there is some linear function $L : \R \rightarrow \R$ such that the rank of
$HF_*(\lambda h(r))$ is at most $L(R)$.
By \cite[Lemma 4.15]{McLean:affinegrowth} and
\cite[Lemma 3.1]{McLean:affinegrowth} we have that the growth rate $\Gamma(M)$
is bounded above by $\varlimsup_x \frac{\log(a(\lambda))}{\log(\lambda)}$ where
$a(\lambda)$ is the rank of $HF_*(\lambda h(r))$.
This implies that $\Gamma(M) \leq \varlimsup_x \frac{\log(L(\lambda))}{\log(\lambda)} \leq 1$.
Hence we have given a bound for $SH_k(M)$ for all $k \notin [1-n,n]$
and also shown that $\Gamma(M) \leq 1$.
\qed

\section{Construction of our exotic Liouville domain} \label{section:exoticliouvilledomain}
In this section we will prove Theorem \ref{theorem:exoticliouvilledomain}.
Here is a statement of this theorem:
{\it In each even dimension greater than $6$ there is a Liouville domain $M$
diffeomorphic to the ball such that $SH_*(M,\Q)$ has infinite rank in each degree}.

From now on our coefficient field will be $\Q$.
We need the following fact:
Let $N$ be any Liouville domain with a choice of trivialization of its canonical bundle
and $l$ any loop in $\partial N$.
We suppose that the dimension of $N$ is greater than $4$.
Then we can attach a Weinstein $2$-handle
along another loop homotopic to $l$ in such a way that the trivialization of the canonical
bundle extends over this handle.

We will not define what a Weinstein handle is here.
The only fact we need to know is that it is a $2$-handle such that the Liouville
domain structure extends over this handle, and also that attaching such a handle
does not change symplectic homology (see \cite{Cieliebak:handleattach}).

\begin{lemma} \label{lemma:increasingdimension}
Let $k$ be any even integer.
Consider the free graded algebra $\K[x,x^{-1},y]$ where $x$ has degree $k$
and $y$ has degree $k+1$.
Let $N$ be a contractible Stein domain of dimension greater than $2$,
then there exists another Stein domain $N'$ such that
\begin{enumerate}
\item $\text{dim}(N') = \text{dim}(N)+2$.
\item $N'$ is contractible.
\item $SH_*(N') = SH_*(N) \otimes \K[x,x^{-1},y]$. 
\end{enumerate}
\end{lemma}
\proof of Lemma \ref{lemma:increasingdimension}.

Notation:
Really by $\C^*$ we mean the Liouville domain associated to $\C^*$
which is the annulus.
Also if we take the product of two Liouville domains $A \times B$, then this is a
manifold with corners. We can smooth the corners to make this a Liouville
domain, but we will just write $A \times B$ for this Liouville domain by abuse of notation.

The set of trivializations of $T\C^*$ is in $1-1$ correspondence with $\Z$.
We normalize so that the trivialization corresponding to $0$
is the natural one coming from viewing $\C^*$ as $\C / \Z$ where we have
the $\Z$ equivariant trivialization of $\C$ induced by the coordinates $a + ib$.
We choose the trivialization corresponding to $\frac{k}{2}$.
From \cite{abbondandoloschwarz:cotangentloopproduct}, we have that
$SH_*(\C^*) = \K[x,x^{-1},y]$.
Really the result in \cite{abbondandoloschwarz:cotangentloopproduct}
uses the trivialization corresponding to $0$ but changing trivialization changes the degree
of $x$ and $y$.
By the statement before this Lemma, we can attach a Weinstein $2$ handle to
$N \times \C^*$ killing the unique generator of $H_1(N \times \C^*) = H_1(\C^*)$ giving us a new Liouville domain $N'$ which is contractible.
We have $SH_*(N') = SH_*(N \times \C^*)$.
By \cite{Oancea:kunneth}, $SH_*(N \times \C^*) = SH_*(N) \otimes SH_*(\C^*)$.
Hence $SH_*(N') = SH_*(N) \otimes \K[x,x^{-1},y]$.
\qed

\proof of Theorem \ref{theorem:exoticliouvilledomain}.
Throughout this proof, our coefficient field $\K$ is equal to $\Q$.
We wish to create a Liouville domain diffeomorphic to the ball
of dimension $2n \geq 8$.
Let $A_1$ be the algebra $\K[x,x^{-1},y]$ where $x$ has degree $0$
and $y$ has degree $1$.
Let $A_2$ be the same algebra but now $x$ has degree $2$ and $y$ has degree $3$.
We start with a contractible Stein domain $D$
of real dimension $4$ with non-trivial symplectic homology
(see \cite[Section 5]{Seidel:biasedview}).
By Lemma \ref{lemma:increasingdimension} there is another contractible Stein domain
$D'$ whose dimension is $\text{dim}(D)+2$ and such that $SH_*(D') = SH_*(D) \otimes A_2$.
This means that $SH_*(D')$ is non-zero in every degree.
We now apply Lemma \ref{lemma:increasingdimension} multiple times to create a contractible Stein domain $M$
of dimension $2n$ whose symplectic homology group is:
$SH_*(D) \otimes A_2 \otimes A_1^{\otimes^{n-3}}$.
Tensoring with $A_1^{\otimes^{n-3}}$ ensures that symplectic homology is now infinitely generated in every degree
as a vector space.
Also $M$ is diffeomorphic to the ball by \cite[Corollary 2.30]{McLean:symhomlef}.
\qed

\bibliography{references}

\end{document}